\documentclass[11pt]{article}
\usepackage{hyperref}
\usepackage{amsfonts,amsmath,amssymb,amsopn}
\usepackage[all]{xy}
\usepackage{vmargin}
\title{The overconvergent site II. Cohomology}
\author{Bernard Le Stum\thanks{bernard.le-stum@univ-rennes1.fr}\\Universit\'e de
Rennes 1}
\date{Version of \today}
\newtheorem{thm}{Theorem}[section]
\newtheorem{prop}[thm] {Proposition}
\newtheorem{cor}[thm] {Corollary}
\newtheorem{lem}[thm] {Lemma}
\newtheorem{dfn}[thm] {Definition}
\def \Id {\mathop{\rm Id}\nolimits}
\def \rig {\mathop{\rm	 rig}\nolimits}
\def \AN {\mathop{\rm AN}\nolimits}
\def \Sch {\mathop{\rm Sch}\nolimits}
\def \Char {\mathop{\rm Char}\nolimits}
\def \an {\mathop{\rm an}\nolimits}
\def \sp {\mathop{\rm sp}\nolimits}

\begin{document}
\maketitle

\begin{center}\textbf{Abstract}\end{center}

We prove that rigid cohomology can be computed as the cohomology of a site analogous to the crystalline site. Berthelot designed rigid cohomology as a common generalization of crystalline and Monsky-Washnitzer cohomology. Unfortunately, unlike the former, the functoriality of the theory is not built-in. We defined somewhere else the ``overconvergent site'' which is functorially attached to an algebraic variety and proved that the category of modules of finite presentation on this ringed site is equivalent to the category of overconvergent isocrystals on the variety. We show here that their  cohomology also coincides.

\tableofcontents

%
%
\section*{Introduction}

Rigid cohomology is a cohomological theory for algebraic varieties over a field $k$ of positive characteristic $p$ with values in vector spaces over a $p$-adic field $K$.
The idea is to embed the given variety $X$ into a proper variety $Y$ and then $Y$ into a smooth formal scheme on the valuation ring $\mathcal V$ of $K$. Then, one considers the limit de Rham cohomology on (strict) neighborhoods of the tube of $X$ inside the generic fibre of $P$.
Classically, one uses Tate's approach to $p$-adic geometry but one may work as well with Berkovich theory and this is what we do here. Anyway, the hard part in the theory is to show that the cohomological spaces obtained this way are independent of the choices (and that they glue when there exists no embedding as above).

There is a relative theory and one may also add coefficients. The coefficients are those of the de Rham theory, namely modules with integrable connexions, on a neighborhood of the tube of $X$, with the extra condition that the connexion must be \emph{overconvergent} (overconvergent isocrystals). Here again, this is a local definition and one must show that it does not depend on the choices. It is also important to remark that glueing is not a satisfactory technic when there is no global embedding as described above. And this is unfortunately the case in general in the relative situation.

In \cite{LeStum07*}, for an algebraic variety $X$ over a formal $\mathcal V$-scheme $S$, we described a ringed site, the \emph{overconvergent site} $(X/S)_{\AN^\dagger}$, and showed that the category of modules of finite presentation on this site is equivalent to the category of overconvergent isocrystals on $X/S$. We prove in this article that cohomology also coincides.
More precisely, we show in Proposition \ref{rigcomp} below that if $p : X \to C$ be a morphism of algebraic varieties over $S$, if $E$ is an overconvergent isocrystal on $X/S$ and $C \hookrightarrow Q$ is a formal embedding over $S$, then
$$
Rp_{\rig} E = Rp_{X/C \subset Q*} E.
$$
where $p_{X/C \subset Q} : (X/S)_{\AN^\dagger} \to (]C[_{Q})_{\an}$ is the canonical morphism.
As a particular case, if $X$ is an algebraic variety over $k$ and $E$ an overconvergent isocrystal on $X/K$, we have for all $i \in \mathbf N$,
$$
H^i_{\rig}(X/K, E) = H^i((X/\mathcal V)_{\AN^\dagger}, E).
$$
Actually, this comparison theorem is a consequence of a slightly more general theorem which states that the cohomology of a module of finite presentation on the overconvergent site can be computed using de Rham resolutions. In order to prove that, we will use linearization of differential operators as in Grothendieck's original article \cite{Grothendieck68}.
If one tries to mimic the classical arguments, one keeps stumbling.
There is no such thing as a projection morphism. The linearization does not give rise to crystals. Linearization is not even local on $X$.
The coefficients do not get out of the de Rham complex.
Some functors are not exact anymore.
Nevertheless, the method works and we get what we want in the end.
The main idea is to use \emph{derived linearization} and the Grothendieck topology.

In section 0, we recall some definitions from (\cite{LeStum06*}). There is nothing new apart from the fact that we will consider the overconvergent site over a general analytic variety $(C, O)$ which is slightly more general that the overconvergent site over a formal scheme $S$ (which corresponds to the case $(S_{k}, S_{K})$).
Also, for technical reasons, we will only consider good analytic varieties.
It should be possible to get rid of this assumption with a systematic use of the Grothendieck topology.

In the first section, we state and prove Berthelot's strong fibration theorem in Berkovich theory.
The main ingredient (invariance of neighborhoods by proper \'etale map) is deduced from general properties of Berkovich analytic varieties.
The rest of the proof follows the original lines.
These results are only used in section 5.

The main point in the second section is a technical base change for finite quasi-immersions that is necessary in order to pull differential operators to the crystal level in the following section.

In section 3, we study differential operators on strict neighborhoods, pull them to differential operators at the crystal level and push them to the analytic site. It is necessary at this point to use derived functors.
Anyway, we show that the cohomology of the derived linearization of a complex with differential operators is captured by the cohomology of the original complex.

Section 4 is totally independent of the previous ones and is devoted to the comparison of the behavior of overconvergence for the analytic and for the Grothendieck topology.
This is necessary because the strong fibration theorem that we proved in section 1 holds only locally for the Grothendieck topology. And we need this result in the proof of the main theorem of the last section.

Section 5 contains the main results. We prove that the cohomology of an overconvergent module of finite presentation can be computed via de Rham cohomology and derive some consequences such as the fact that this cohomology coincides with rigid cohomology.

In Appendix A, we show that, even if the usual topology of an analytic variety is not rich enough to recover the Zariski topology of its reduction, the cohomology of crystals on the overconvergent site is local for the Zariski topology of $X$.

Finally, in Appendix B, we explain how to compute the cohomology of general complexes of abelian sheaves using \v{C}ech-Alexander techniques.
Here again, it is necessary to work in the derived category and define the \emph{derived \v{C}ech-Alexander complex}.

%
%

%
%
\section*{Acknowledgments}

I want to thank here Pierre Berthelot: working close to him is an incredible chance. I dropped the linearization method for one year, looking for other ways to prove the comparison theorem, and he convinced me to try again. 
Many thanks also to Antoine Ducros who is so fluent in Berkovich theory. 
I bothered him several times with elementary questions (elementary for him !).

%
%
\section*{Conventions}

Throughout this paper, $K$ is a \emph{non trivial} complete ultrametric field with valuation ring $\mathcal V$, maximal ideal $\mathfrak m$ and residue field $k$.
When using Berkovich theory, we consider only \emph{strictly} analytic varieties.
I did not think at all about the trivial valuation case.
Finally, complexes are always assumed to be bounded below.

%
%

%
%

\section{The overconvergent site}

For the convenience of the reader, we recall in this section some results and definitions from \cite{LeStum06*}.

\subsection{Formal schemes}

When the valuation is \emph{discrete}, a \emph{formal $\mathcal V$-scheme} is a locally finite union of formal affine $\mathcal V$-schemes of finite type.
In general, we restrict to formal affine schemes defined by an ideal of the form $\mathfrak a\mathcal V\{T_{1}, \ldots, T_{n}\} + I$ with $I$ finitely generated.
The formal scheme is said \emph{admissible} if it is $\mathcal V$-flat.
We may (and will) see the category of algebraic $k$-varieties as a full subcategory of the category of formal $\mathcal V$-schemes.

We will systematically consider the \emph{special fibre} $P_{k}$ of a formal $\mathcal V$-scheme $P$, which is an algebraic $k$-variety, as well as its \emph{generic fibre} $P_{K}$, which is an analytic variety (in the sense of Berkovich).
There is an immersion $P_{k} \hookrightarrow P$ which is a homeomorphism and a  \emph{specialization map} $\sp : P_{K} \to P$ (which is anticontinuous).

A \emph{formal embedding} is a locally closed immersion $X \hookrightarrow P$ of an algebraic $k$-variety into a formal $\mathcal V$-scheme.
The tube of $X$ in $P$ is
$$
]X[_{P} := \{x \in P, \quad \sp(x) \in X\}.
$$
The formal embedding is said \emph{good} if any point of $]X[_{P}$ has an affinoid neighborhood in $P_{K}$.
A morphism of formal embeddings is simply a pair of compatible morphisms.

\subsection{The big overconvergent site}
The \emph{(big) overconvergent site} $\AN^\dagger(\mathcal V)$ is defined as follows.

An \emph{object} is a ``couple'' $(X, V)$ made of a good formal embedding $X \hookrightarrow P$ of a $k$-variety into a $\mathcal V$-formal scheme and a morphism $\lambda : V \to P_K$ from an analytic variety to the generic fibre of $P$.
We call such a couple an \emph{analytic variety} over $\mathcal V$.
It is convenient to consider the corresponding tube
$$
]X[_{V} = \lambda^{-1}(]X[_{P})
$$
and define a \emph{(strict) neighborhood} of $X$ in $V$ as an analytic domain which is a neighborhood of $]X[_{V}$ in $V$.

We call an analytic variety $(X, V)$ over $\mathcal V$, \emph{good} if any point of $]X[_{V}$ has an affinoid neighborhood in $P_{K}$.
Unlike in \cite{LeStum06*}, we will consider here \emph{only} good analytic varieties over $\mathcal V$.
The main reason is that we need coherent cohomology to coincide wether we use the analytic topology or the Grothendieck topology. However, most results will still be valid without this assumption.

A \emph{morphism} $(X', V') \to (X, V)$ of analytic varieties over $\mathcal V$, is a couple of ``geometrically pointwise compatible'' morphisms
$$
f : X' \to X, \quad u : V'' \to V
$$
where $V''$ is a strict neighborhood of $X'$ in $V'$.
We need to make more precise the expression ``geometrically pointwise compatible'' : it means that 
$$
\forall x \in ]X[_{{V'}}, \quad \sp(\lambda(u(x)) = f(\sp(\lambda'(x))
$$
and this still holds after any isometric extension $K \hookrightarrow K'$.
The reader who does not feel comfortable with our definition of morphisms above should note that any commutative diagram like this
$$
\xymatrix{X' \ar@{^{(}->}[r] \ar[d]_f & P' \ar[d]^{v} & P'_{K} \ar[l] \ar[d]_{v_{K}} &V' \ar[l] \ar[d]^u \\ X \ar@{^{(}->}[r] & P & P_{K} \ar[l] &V \ar[l] }
$$
gives rise to a morphism $(X', V') \to (X, V)$.
Moreover, up to an isomorphism upstairs, any morphism of analytic varieties $(X', V') \to (X, V)$ over $\mathcal V$ has this form.
Note also that, in general, a morphism of analytic varieties
$$
(f, u) : (X', V') \to (X, V)
$$
induces a morphism on the tubes
$$
]f[_u : ]X'[_V' \to ]X[_V.
$$
When there is no risk of confusion, we will simply write $u$.
When $u$ is the identity, we will write $]f[$.

Finally, a \emph{covering} of an analytic variety $(X, V)$ over $\mathcal V$ is simply a family $(X, V_{i})$ where $V' := \cup V_{i}$ is an open covering of a neighborhood of $X$ in $V$.

\subsection{Restricted sites}

 We will systematically make use of \emph{restriction} (it is usually called localization) in the following sense.
If $C$ is a category and $T$ is a presheaf on $C$, then $C_{/T}$ is the category of pairs $(X, s)$ with $X \in C$ and $s \in T(X)$.
Morphisms are the obvious ones: morphisms in $C$ that are compatible with the given sections.
This applies in particular to representable presheaves in which case, we obtain the usual notion of restricted category (all the objects above another one  with compatible morphisms).

There is always a forgetful functor $j_{T} : C_{/T} \to C$ which defines a morphism of toposes $j_{T} : \widetilde{C_{/T}} \to \widetilde{C}$ (with three adjoint functors $j_{!}, j^{-1}, j_{*}$) when $C$ is a site.
It is important to remark that, on abelian sheaves, there exists also a left adjoint $j^{ab}_{!}$ which is exact so that $j^{-1}$ preserves injectives.
Finally, if $\widetilde{C}$ denotes the topos associated to $C$ and $\widetilde {T}$ the sheaf associated to $T$, we have $\widetilde{C_{/T}} \simeq \widetilde {C}_{/\widetilde{T}}$.

Any morphism of presheaves $u : T' \to T$ induces a functor $u_{C} : C_{/T'} \to C_{/T}$ which gives a morphism of toposes when $C$ is a site.
Actually, this is a particular case of a restriction map as before because we may see the morphism $u$ as a presheaf on $C_{/T}$.
Finally, as for open embeddings in topological spaces, if we are given a cartesian diagram of presheaves
$$
\xymatrix{T_{1} \times_{T} T_{2} \ar[r]^-{p_{1}} \ar[d]^{p_{2}}  & T_{1} \ar[d]^{u_{1}} \\ T_{2} \ar[r]^{u_{2}} & T,}
$$
then we have a general base change theorem that reads
$$
u_{1C}^{-1}u_{2C*} \mathcal F \simeq p_{1C*}p_{2C}^{-1} \mathcal F.
$$
For complexes of abelian sheaves, we also have
$$
u_{1C}^{-1}Ru_{2C*} \mathcal F \simeq Rp_{1C*}p_{2C}^{-1} \mathcal F
$$
since inverse images preserve injectives.

\subsection{Functoriality}

If $T$ is an \emph{analytic presheaf} on $\mathcal V$, which means a presheaf on $\AN^\dagger(\mathcal V)$, we will write
$$
\AN^\dagger(T) := \AN^\dagger(\mathcal V)_{/T}
$$
for the restricted site and $T_{\AN^\dagger}$ for the corresponding topos.
An object of $\AN^\dagger(T)$ will be called an analytic variety over $T$.
We will apply this, for example, to the case of a representable presheaf $(X, V)$.
We will write $(X \subset P)$ instead of $(X, P_{K})$ when $(X \subset P)$ is a good formal embedding and $S$ instead of $(S_{k} \subset S)$ when $S$ is a formal $\mathcal V$-scheme.
Thus, we get the sites $\AN^{\dagger}(X, V)$, $\AN^{\dagger}(X \subset P)$ or $\AN^{\dagger}(S)$ and the corresponding topos.

Any morphism of analytic presheaves $u : T' \to T$ induces by restriction a morphism of toposes
$$
u_{\AN^\dagger} : T'_{\AN^\dagger} \to T_{\AN^\dagger}.
$$ 
In particular, if $(X, V)$ is an analytic variety over an analytic presheaf $T$, there is a restriction map
$$
j_{X,V} : (X,V)_{\AN^\dagger} \to T_{\AN^\dagger}
$$
and if $(f, u) : (X', V') \to (X, V)$ is a morphism, there is also a restriction map
$$
u_{\AN^\dagger} : (X',V')_{\AN^\dagger} \to (X,V)_{\AN^\dagger}.
$$

There is also a morphism of sites which is very important
$$
\xymatrix@R0cm{\AN^\dagger(X, V) \ar[rr]^{\varphi_{X,V}} && \mathrm{Open}(]X[_{V})\\Ê(X, V') && \ar[ll] ]X[_{V'}.}
$$
As morphism of topos, $\varphi_{X,V}$ has a section $\psi_{X,V}$ with $\varphi_{X,V*} = \psi_{X,V}^{-1}$.
When $T$ is an analytic presheaf on $\mathcal V$, $\mathcal F$ an analytic sheaf on $T$, and $(X, V)$ an analytic variety on $T$, we will write
$$
\mathcal F_{X,V} := \psi_{X,V}^{-1}j_{X,V}^{-1} \mathcal F
$$
and call it the \emph{realization} of $\mathcal F$ on $(X, V)$.
Note that this realization functor $\mathcal F \mapsto \mathcal F_{X,V}$ is exact and preserves injectives.

\subsection{The relative overconvergent site}

We will define now the analytic presheaf ``$X/(C,O)$'' as well as the corresponding site and topos.
This is a generalization of the situation we considered in \cite{LeStum06*}, that corresponds here to the case $(C , O) = (S_{k}, S_{K}) =: S$ for a formal $\mathcal V$-scheme $S$.

If $\Sch/k$ denotes the category of algebraic varieties over $k$ with the coarse topology, there is an obvious morphism of sites (induced by the forgetful functor)
$$
\xymatrix@R0cm{
 \Sch/k \ar[rr]^I && \AN^{\dagger}(\mathcal V) \\ X && \ar[ll] (X, V).
}
$$
If $(C, O)$ is a an analytic variety over $\mathcal V$, then $I$ extends formally to a morphism of sites
$$
I_{C,O} : \Sch/C \to  \AN^{\dagger}(C, O).
$$
We may also consider the restriction morphism
$$
j_{C,O} : (C, O)_{\AN^{\dagger}} \to \mathcal V_{\AN^{\dagger}}
$$
and we define, when  $X$ is any algebraic variety over $C$, the analytic presheaf
$$
X/(C,O) := j_{C,O\,!}I_{C,O\,*} X.
$$
An object is just a pair of morphisms
$$
(U \to X, \quad V \to O)
$$
such that the composite
$$
(U \to X \to C, \quad V \to O)
$$
is a morphism of analytic varieties.
And morphisms $(U', V') \to (U, V)$ are simply usual morphisms of analytic varieties  over $(C, O)$ where $U' \to U$ is actually a morphism over $X$. 

Note that if $(C, O)$ is an analytic variety over $S$, we have an isomorphism of presheaves
$$
X/(C,O) \simeq X/S \times_{C/S} (C, O).
$$

\subsection{Cohomology}

If $(C, O)$ is an analytic variety over $\mathcal V$ and $p : X \to C$ a morphism of algebraic varieties, we may consider the projection
$$
\xymatrix@R0cm{p_{X/C, O} : & (X/C, O)_{\AN^{\dagger}} \ar[r] & (C, O)_{\AN^\dagger} \ar[r] & (]C[_{O})_{\an}. \\
&X/(C, O') && O' \ar@{|->}[ll]}
$$
Deriving this functor (for abelian sheaves) gives \emph{absolute cohomology}.

On the other hand, any morphism $f : X' \to X$ of algebraic varieties over $C$, induces a morphism of analytic presheaves $f : X'/(C,O) \to X/(C,O)$ which gives a morphism of toposes
$$
f_{\AN^\dagger} : (X'/C, O)_{\AN^\dagger} \to (X/C,O)_{\AN^\dagger}.
$$
Deriving this functor (for abelian sheaves) gives \emph{relative cohomology}.

Base change for restriction maps implies that if $\mathcal F'$ is a complex of abelian sheaves on $X'/(C,O)$ and $(U, V)$ an analytic variety over $X/(C,O)$, then
$$
(R f_{\AN^\dagger*} \mathcal F')_{U, V} = Rp_{X' \times_{X} U/(U, V)*} \mathcal F'_{|X' \times_{X} U/(U, V)}.
$$
This shows that relative cohomology can be recovered from absolute cohomology.
Of course, conversely, we have for a complex of abelian sheaves $\mathcal F$ on $X/(C, O)$,
$$
Rp_{X/C, O*} \mathcal F = (R p_{\AN^\dagger*} \mathcal F)_{C,O}
$$
if $p : X \to C$ denotes the structural map as above.
Therefore, absolute cohomology can also be derived from relative cohomology.

\subsection{Crystals}

If $T$ is an analytic presheaf on $\mathcal V$, we will consider the sheaf of \emph{overconvergent functions} on $T$ defined by
$$
\mathcal O^\dagger_{T}(X, V) = \Gamma(]X[_{V}, i_{X}^{-1}\mathcal O_{V})
$$
where $i_{X} : ]X[_{V} \hookrightarrow V$ denotes the embedding.
An \emph{overconvergent module} on $T$ will be an $\mathcal O^\dagger_{T}$-module $E$ and it will be called a \emph{crystal} if the transition maps are bijective
$$
u^\dagger E_{X,V} := i_{X'}^{-1}u^*i_{X*}E_{X,V} \simeq E_{X', V'}
$$
whenever $(f, u) : (X', V') \to (X, V)$ is a morphism in $\AN^\dagger(T)$.

Restriction morphisms are automatically morphisms of ringed topos: if $u : T' \to T$ is any morphism of analytic presheaves, we have
$$
(u_{C}^{-1}, u_{C*}) :  (T'_{\AN^\dagger}, \mathcal O^\dagger_{T'}) \to (T_{\AN^\dagger}, \mathcal O^\dagger_{T}).
$$
In particular, inverse image is exact and preserves injectives.

Realization becomes also a morphism of ringed sites (see proposition 7.4 of \cite{LeStum06*})
$$
(\varphi_{X,V}^\dagger, \varphi_{X,V*}) : (\textrm{AN}^\dagger(X,V), \mathcal O_{X,V}^\dagger) \to (]X[_V, i_X^{-1}\mathcal O_V).
$$
It induces an equivalence between overconvergent crystals on $(X, V)$ and $i_X^{-1}\mathcal O_V$-modules.
Moreover, finitely presented overconvergent modules on $(X, V)$ correspond to coherent $i_X^{-1}\mathcal O_V$-modules.
Finally, note that the functor $\varphi_{X,V*}$ is exact but that $\varphi_{X,V}^\dagger$ is not left exact in general.

%
%
\section{The strong fibration theorem}

We reprove here in the context of Berkovich theory the strong fibration theorem of Berthelot.
The strategy of the proof is exactly the same as Berthelot's.

The following result is the key point.
The classical proof is quite involved (Theorem 1.3.5 of \cite{Berthelot96c*} or Theorem 3.4.12 \cite{LeStum07*}).
We see it here as an add on to section 2 of \cite{LeStum06*}.

\begin{prop}\label{isom}
If $u : P' \to P$ be a morphism of good admissible embeddings of an algebraic variety $X$ that is proper and \'etale at $X$, it induces an isomorphism of analytic varieties $(X, P_{K}') \simeq (X, P_{K})$.
\end{prop}

\textbf {Proof : }
It follows from corollary 2.6 of \cite{LeStum06*} that $u$ induces an \'etale morphism $V' \to V$ between neighborhoods of $X$ in $P$ and $P'$ respectively. Moreover, it follows from lemma 44 of \cite{Berkovich99} that $u$ induces an isomorphism $]X[_{P'} \simeq ]X[_{P}$ on the tubes.
Now, the result follows from proposition 4.3.4 of \cite{Berkovich93}.
$\Box$

If $X \subset P$ is a formal embedding, we may embed $X$ in $\widehat{\mathbf A}^n_{P}$ using the zero section in order to get another  formal embedding with a canonical projection
$$
(X \subset \widehat{\mathbf A}^n_{P}) \to (X \subset P).
$$
This morphism is proper and smooth at $X$.
The strong fibration theorem says that, locally, any morphism that is proper and smooth at $X$ looks like this one.
More precisely, we have:

\begin{thm}[Strong Fibration Theorem] \label{strong}
Let $u : P' \to P$ be a morphism of good admissible embeddings of an algebraic variety $X$ that is proper and smooth at $X$.
Then, locally for the Zariski topology on $X$ and for the Grothendieck topology on $P_{K}$, the corresponding analytic variety $(X, P'_{K})$ over $(X, P_{K})$ is isomorphic to $(X, \widehat{\mathbf A}^n_{P_{K}})$.
\end{thm}

Almost all the ideas in the following proof are already in Berthelot's original preprint \cite{Berthelot96c*}.

\textbf {Proof : }
Using corollary 3.5 of \cite{LeStum06*}, we may assume that the locus at infinity of $X$, both in $P$ and in $P'$, is a the support of a divisor.
Also, being local for the Grothendieck topology of $P_{K}$, the question is local for the Zariski topology of $P$.
In particular, we may assume that $P$ is affine and that the locus at infinity of $X$ in $P$ is actually a hypersurface.
Finally, removing extra components, we may also assume that $P'$ is quasi-compact.

Now, let $Y$ be the Zariski closure of $X$ in $P'_{k}$.
By hypothesis, the map $Y \to P_{k}$ induced by $u$, is proper.
Using corollary 3.4 of \cite{LeStum06*}, we may therefore assume thanks to Chow's lemma (corollary I.5.7.14 of \cite{GrusonRaynaud71}) that the induced map $Y \to P_{k}$ is projective.
Thus, there exists a commutative diagram
$$
\xymatrix{& Y \ar@{^{(}->}[r] \ar[dd] & \mathbf P^N_{P} \ar[dd]^p \\ X \ar@{^{(}->}[ur] \ar@{^{(}->}[dr] && \\ & P_{k} \ar@{^{(}->}[r] & P.}
$$
The closed immersion $X \hookrightarrow p^{-1}(X) = \mathbf P_{X}^N$ is a section of the canonical projection which is smooth.
Since the question is local on $X$, we may assume that there exists an open subset $U$ of $ \mathbf P_{X}^N$ such that $X$ is defined in $U$ by a regular sequence.
If $\overline X$ denotes the Zariski closure fo $X$ in $P$, we may assume that $U = D^+(\bar s) \cap \mathbf P_{\overline{X}}^N$ with $s \in \Gamma(\mathbf P_{P}^N, \mathcal O(m))$ for some $m$ and that the regular sequence is induced by $t_{1}, \dots , t_{d} \in \Gamma(\mathbf P_{P}^N, \mathcal O(n))$ for some $n$.
Then, by construction, the induced morphism
$$
v: P'' := V^+(t_{1}, \cdots, t_{d}) \to P
$$
is proper and \'etale at $X$ and lifts the map $Y \to \overline X$.
We embed $Y$ in the product $P''' := P'' \times_{P} P'$ and consider the following cartesian diagram of embeddings of $X$:
$$
\xymatrix{P' \ar[d]_u & P''' \ar[d]^{p_{2}} \ar[l]_{p_{1}}Ê\\ P & P'' \ar[l]^v.}
$$
Proposition \ref{isom} tells us that both horizontal arrows induce an isomorphism on the corresponding analytic varieties.
The question being local on $P_{K}$, we may therefore replace $u$ by $p_{2}$ and assume that $u$ induces an isomorphism on $Y \simeq \overline X$.
The question is henceforth local on $Y$ too. Good.

Since the question is local on $X$, we may assume that the conormal sheaf $\check{\mathcal N}$ of $X$ in $u^{-1}(X)$ is free.
And we can replace $P'$ by an affine neighborhood of $Y$.
We may then lift a basis of $\check{\mathcal N}$ to a sequence $s_{1}, \ldots, s_{n}$ of elements of the ideal of $Y$ in $P'$.
By construction, they induce a morphism $P'_{K} \to \widehat{\mathbf A}^n_{P_{K}}$  of good admissible embeddings of $X$ which is proper \'etale at $X$.
And we get as expected an isomorphism of analytic varieties over $(X, P_{K})$ using again proposition \ref{isom}.
$\Box$

In order to apply this theorem, we will have to localize on $X$ when $(X, V)$ is an analytic variety on $\mathcal V$.
This is possible thanks to the next result.

\begin{prop} \label {Zar}
Let $(X, V)$ be an analytic variety over $\mathcal V$ and $\{X_{k}\}_{k \in I}$ a locally finite open covering of $X$ with inclusion maps
$$
\alpha_{k} : X_{k} \hookrightarrow X, \quad \alpha_{kl} : X_{k} \cap X_{l} \hookrightarrow X, \quad \ldots
$$
If $\mathcal F$ is an abelian sheaf on $]X[_{V}$, there is a long exact sequence
$$
0 \to \mathcal F \to \prod_{k} ]\alpha_{k}[_{V*}]\alpha_{k}[_{V}^{-1}  \mathcal F \to  \prod_{k,l} ]\alpha_{kl}[_{V*}]\alpha_{kl}[_{V}^{-1}  \mathcal F \to \cdots 
$$
\end{prop}

\textbf {Proof : }
This assertion is local on $V$.
Since the specialization map is anticontinuous, the tubes of the irreducible components of $X$ form an open covering of $]X[_{V}$. We may thus assume that $X$ is irreducible and therefore quasi-compact, in which case our covering is finite. Again, since the specialization map is anticontinuous, we obtain a finite closed covering and our sequence is simply the Mayer-Vietoris sequence for this closed covering.
$\Box$

\begin{cor} \label{locx}
Let  $(f, u) : (X', V') \to (X, V)$ be a morphism of analytic varieties, $\mathcal F'$ a complex of abelian sheaves on $]X'[_{V'}$ and $\{X'_{k}\}_{k \in I}$ a locally finite open covering of $X'$.
Then, there is a spectral sequence
$$
E_{1}^{r,s} = \bigoplus R^s u_{|]X_{\underline k}'[_{V'}*} \mathcal F'_{|]X_{\underline k}'[_{V'}*}\Rightarrow R^{r+s}u_{*}\mathcal F'.
$$
\end{cor}

\textbf {Proof : }
Of course, we use the fact that all maps $]\alpha[_{\underline k} : ]X'[_{\underline k}[_{V'} \to ]X'[_{V'}$ are inclusions of closed subsets and therefore that all $]\alpha[_{\underline k*}$ are exact.
$\Box$

%
%
\section{Functoriality}

In this section, we prove some functoriality results that are used in the next one.

Recall that when $(f, u) : (X', V') \to (X, V)$ is a morphism of analytic varieties over $\mathcal V$, there is a restriction map
$$
u_{\AN^\dagger} : (X', V')_{\AN^\dagger} \to (X, V)_{\AN^\dagger}
$$
which is naturally a morphism of ringed toposes.

\begin{prop} \label{restr}
Let 
$$
\xymatrix{(Y', W') \ar[rr]^{(f', u')}Ê\ar[d]_{(g', v')}&& (Y, W) \ar[d]^{(g, v)} \\ (X', V') \ar[rr]^{(f, u)} && (X, V)}
$$
be a cartesian diagram of analytic varieties over $\mathcal V$.
If $\mathcal F'$ is an analytic sheaf on $(X', V')$, we have
$$
(u_{\AN^\dagger*} \mathcal F')_{Y, W} = u'_{*}\mathcal F_{Y', W'}.
$$
Actually, if $\mathcal F'$ is a complex of abelian sheaves, we even get
$$
(Ru_{\AN^\dagger*} \mathcal F')_{Y, W} = Ru'_{*}\mathcal F'_{Y', W'}.
$$
\end{prop}

\textbf {Proof : }
The first assertion follows from the fact that $u_{\AN^\dagger}^{-1}(Y, W) = (Y' , W')$ and therefore
$$
\Gamma((Y, W), u_{\AN^\dagger*} \mathcal F') = \Gamma(u_{\AN^\dagger}^{-1}(Y, W), \mathcal F') = \Gamma((Y' , W'), \mathcal F') $$
$$
= \Gamma(]Y'[_{W'}, \mathcal F'_{Y', W'}) =  \Gamma(]X'[_{V'}, u'_{*}\mathcal F'_{Y', W'}).
$$
For the second one, it is sufficient to recall that the realization functors are exact and preserve injectives.
$\Box$

Recall that, if $(X, V)$ is an analytic variety over $\mathcal V$, there is a morphism of ringed sites
$$
(\varphi_{X,V}^\dagger, \varphi_{X,V*}) : (\textrm{AN}^\dagger(X,V), \mathcal O_{X,V}^\dagger) \to (]X[_V, i_X^{-1}\mathcal O_V).
$$
Also, any morphism $(f, u) : (X', V') \to (X, V)$ of analytic varieties over $\mathcal V$ induces a morphism of ringed spaces
$$
(u^\dagger, u_{*}) : (]X'[_{V'}, i_{X',V'}^{-1} \mathcal O_{V'}) \to (]X[_V, i_{X,V}^{-1} \mathcal O_V).
$$
When $\mathcal F$ is an $ i_{X,V}^{-1} \mathcal O_V$-module, we have
$$
(\varphi_{X,V}^{\dagger} \mathcal F)_{X',V'} = u^\dagger \mathcal F.
$$

\begin{prop} \label{comdag}
If $(f, u) : (X', V') \to (X, V)$ is a morphism of analytic varieties over $\mathcal V$ and $\mathcal F$ is an $ i_{X,V}^{-1} \mathcal O_V$-module, there is a canonical isomorphism
$$
u_{\AN^\dagger}^{-1}\varphi_{X,V}^{\dagger} \mathcal F = \varphi_{X',V'}^{\dagger} u^\dagger \mathcal F.
$$
\end{prop}

\textbf {Proof : }
There is a commutative diagram of ringed toposes
$$
\xymatrix{(]X'[_{V'})_{\an} \ar[d]_u && (X', V')_{\AN^\dagger} \ar[ll]_{\varphi_{X',V'}} \ar[d]^{u_{\AN^\dagger}} & \\(]X[_{V})_{\an} && (X, V)_{\AN^\dagger}. \ar[ll]^{\varphi_{X,V}} & \Box}
$$

Recall from \cite{Berkovich93}, Definition 4.3.3, that a \emph{finite quasi-immersion} of analytic varieties is a finite morphism which is a homeomorphism on its image with purely inseparable residue field extensions. This property is stable under composition and base change.

\begin{dfn} \label{findef}
A morphism $(f, u) : (X', V') \to (X, V)$ of analytic varieties  over $\mathcal V$ is \emph{finite} (resp. a \emph{finite quasi-immersion}) if, after replacing $V$ and $V'$ by strict neighborhoods of $X$ and $X'$ respectively, we have
\begin{enumerate}
\item $u$ is finite (resp. a finite quasi-immersion)
\item $u^{-1}(]X[_{V}) = ]X'[_{V'}$
\end{enumerate}
\end{dfn}

\begin{prop} \label{finm}
Let $(f, u) : (X', V') \to (X, V)$ be a finite quasi-immersion of analytic varieties  over $\mathcal V$ and $\mathcal F'$ an $i_{X'}^{-1}\mathcal O_{V'}$-module.
Then the adjunction map
$$
\varphi_{X,V}^{\dagger}u_{*}\mathcal F' \to u_{\AN^\dagger*}\varphi_{X',V'}^{\dagger}\mathcal F'
$$
is an isomorphism. Moreover, we have
$$
R^qu_{*}\mathcal F' = 0 \quad \mathrm{and} \quad R^qu_{\AN^\dagger*}\varphi_{X',V'}^{\dagger}\mathcal F' = 0 \quad \mathrm{for} \quad q > 0.
$$
\end{prop}

There should be an analogous result for finite morphisms and coherent sheaves but we will not need it.

\textbf {Proof : }
Of course, we choose $V$ and $V'$ such that $(f, u)$ satisfies the conditions of Definition \ref{findef}.
Let $(g, v) : (Y, W) \to (X, V)$ be any morphism of analytic varieties, $(g', v') : (Y', W') \to (X', V')$ the pull-back along $(f, u)$, and $(f', u') : (Y', W') \to (Y, W)$ the corresponding map.
We are in the situation of Proposition \ref{restr} and it follow that
$$
(Ru_{\AN^\dagger*} \varphi_{X',V'}^{\dagger}\mathcal F')_{Y,W} = Ru'_{*}(\varphi_{X',V'}^{\dagger}\mathcal F')_{Y', W'}.
$$
In particular, in order to prove the last assertion, it is actually sufficient to show that
$$
R^qu_{*}\mathcal F' = 0 \quad \mathrm{for} \quad q > 0
$$
and then apply it to $u'_{*}$ and $(\varphi_{X',V'}^{\dagger}\mathcal F')_{Y', W'}$ (finite quasi-immersions are preserved by base change).

We also have to show that we always have
$$
(\varphi_{X,V}^{\dagger}u_{*}\mathcal F')(Y,W) \simeq (u'_{*}(\varphi_{X',V'}^{\dagger}\mathcal F')_{(Y', W')})(Y,W).
$$
And this means that
$$
\Gamma(]Y[_{W}, v^\dagger u_{*}\mathcal F') \simeq \Gamma(]Y[_{W}, u'_{*} {v'}^\dagger \mathcal F').
$$
We are therefore reduced to checking that
$$
v^\dagger u_{*}\mathcal F' \simeq u'_{*} {v'}^\dagger \mathcal F'
$$
on $]Y[_{W}$.

Of course, it is sufficient to consider sheaves of the form $i_{X'}^{-1}\mathcal F'$.
Thus, changing notations, we want to show that, if $\mathcal F'$ is an $\mathcal O_{V'}$-module, we have
$$
v^\dagger u_{*}i_{X'}^{-1}\mathcal F' \simeq u'_{*} {v'}^\dagger i_{X'}^{-1} \mathcal F'
$$
and
$$
R^qu_{*} i_{X'}^{-1}  \mathcal F' = 0 \quad \mathrm{for} \quad q > 0.
$$
Since $u$ is just the inclusion of a closed subset, it follows from the second condition of Definition \ref{findef} that the cohomology vanishes and that $u_{*}i_{X'}^{-1}\mathcal F' = i_{X}^{-1}u_{*}\mathcal F'$.
By definition, we obtain
$$
v^\dagger u_{*}i_{X'}^{-1}\mathcal F' = v^\dagger i_{X}^{-1}u_{*}\mathcal F' = i_{Y}^{-1} v^* u_{*}\mathcal F'.
$$
For the same reasons (definition and compactness of $u'$) we first get that ${v'}^\dagger i_{X'}^{-1} \mathcal F' = i_{Y'}^{-1}  {v'}^* \mathcal F'$ and then
$$
u'_{*}{v'}^\dagger i_{X'}^{-1} \mathcal F' = u'_{*} i_{Y'}^{-1}  {v'}^* \mathcal F' = i_{Y}^{-1}u'_{*}{v'}^* \mathcal F'.
$$
It is therefore sufficient to show that
$$
v^* u_{*}\mathcal F' = u'_{*}{v'}^* \mathcal F'
$$
on $W$.
This can be checked on each fiber. If $x \in W$ is not in the image of $W'$, then $v(x)$ is not in the image of $V'$ and both sides vanish.
If $x \in W'$, then
$$
(v^* u_{*}\mathcal F')_{u'(x)} = \mathcal O_{W,u'(x)} \otimes_{\mathcal O_{V,vu'(x)}} (u_{*}\mathcal F')_{vu'(x)} = 
$$
$$
= \mathcal O_{W,u'(x)} \otimes_{\mathcal O_{V,uv'(x)}} (u_{*}\mathcal F')_{uv'(x)} = \mathcal O_{W,u'(x)} \otimes_{\mathcal O_{V,uv'(x)}} \mathcal F'_{v'(x)} 
$$
At last, we use the fact that $u$ is finite to obtain
$$
u'_{*}{v'}^* \mathcal F'_{u'(x)} = ({v'}^* \mathcal F')_{x} = \mathcal O_{W',x} \otimes_{\mathcal O_{V',v'(x)}} \mathcal F'_{x}
$$
$$
= (\mathcal O_{W,u'(x)} \otimes_{\mathcal O_{V,vu'(x)}} \mathcal O_{V',v'(x)}) \otimes_{\mathcal O_{V',v'(x)}}\mathcal F'_{x} = \mathcal O_{W,u'(x)} \otimes_{\mathcal O_{V,vu'(x)}} \mathcal F'_{v'(x)}. \quad \Box
$$

%
%
\section{Differential operators}

In this section, we fix a morphism $(f, u) : (X, V) \to (C, O)$ of analytic varieties over $\mathcal V$.

We will consider the infinitesimal neighborhood of oder $n$ of $V$ which is the analytic subvariety $V^{(n)}$ defined by $\mathcal I^{n+1}$ if $V$ is defined by $\mathcal I$ in $V \times_{O} V$.
The diagonal embedding $\delta :  V \hookrightarrow V \times_{O} V$ and the projections $p_{i} : V \times_{O} V \to V$ induce morphisms $\delta^{(n)} :  V \hookrightarrow V^{(n)}$ and $p^{(n)}_{i} : V^{(n)} \to V$.
From $p_{i} \circ \delta = \Id_{V}$, we deduce that $p^{(n)}_{1}$ and $p^{(n)}_{2}$ are identical as continuous maps (both are inverse to the homeomorphism $\delta^{(n)}$). 
Actually, if we use $\delta^{(n)}$ to identify the underlying topological spaces of $V^{(n)}$ and $V$, then $p_{i}^{(n)}$ becomes the identity as continuous maps both for $i = 1$ and $i = 2$ so that $p_{2*}^{(n)} = p_{1*}^{(n)}$.
Of course, the projections act differently on sections so that $p^{(n)*}_{2} \neqÊp^{(n)*}_{1}$.

We will embed $X$ in $V^{(n)}$ using $\delta^{(n)}$ and consider the resulting analytic variety $(X, V^{(n)})$ over $(C, O)$ with the projection maps $p^{(n)}_{i} : (X, V^{(n)}) \to (X, V)$.
Again these maps do nothing on the underlying topological spaces and only play a role on sections.

\begin{dfn}  An $i_C^{-1}\mathcal O_{O}$-linear morphism $d : \mathcal F \to \mathcal G$ between two $i_X^{-1}\mathcal O_{V}$-modules is a \emph{differential operator (of order at most $n$)} if it factors as
$$
\xymatrix{\mathcal F \ar[r]^-{p^{(n)*}_{2}} & p^{(n)}_{2*}p^{(n)\dagger}_{2} \mathcal F \ar@{=}[r] & p^{(n)}_{1*}p^{(n)\dagger}_{2} \mathcal F \ar[r]^-{\overline d} & \mathcal G.}
$$
where $\overline d$ is $i_X^{-1}\mathcal O_{V}$-linear
\end{dfn}

Note that the equality sign is \emph{not} linear in general.
Note also that $\overline d$ is unique: locally, we have
$$
\bar d (\overline {f \otimes g} \otimes m) = fd(gm).
$$
The composition of two differential operators is a differential operator.
More precisely, if $d : \mathcal F \to \mathcal G$ and $d' : \mathcal G \to \mathcal H$ are two differential operators of order at most $n$ and $m$ respectively, then $d' \circ d$ is a differential operator of order at most $m+n$ and
$$
\overline{d' \circ d} = \overline{d'} \circ p^{(m)\dagger}_{2}(\overline d).
$$

We will mainly be concerned with differential operators of order at most $1$:
it just means that, locally, we have
$$
d(fgm) + fgd(m) = fd(gm) + gd(fm).
$$
The standard example of differential operators of order at most $1$ is given by the differentials in a de Rham complex
$$
\mathcal F \otimes_{i_{X}^{-1}\mathcal O_{V} } i_{X}^{-1}\Omega^\bullet_{V/O}
$$
when $\mathcal F$ is an $i_{X}^{-1}\mathcal O_{V} $-module with an integrable connexion.

We can extend this notion of differential operator to overconvergent modules as follows:

\begin{dfn}  An $\mathcal O_{C,O}^\dagger$-linear morphism $d : E \to E'$ between two $\mathcal O_{X,V}^\dagger$-modules is a \emph{differential operator (of order at most $n$)} if it factors as the composite of two $\mathcal O_{X,V}^\dagger$-linear morphisms
$$
\xymatrix{E \ar[r]^-{p^{(n)*}_{2\AN^\dagger}} & p^{(n)}_{2\AN^\dagger*}p^{(n)^{-1}}_{2\AN^\dagger} E \ar@{=}[r] & p^{(n)}_{1\AN^\dagger*}p^{(n)^{-1}}_{2\AN^\dagger} E \ar[r]^-{\overline d} & E'.}
$$
\end{dfn}

Of course, here again, the equality sign in this definition is \emph{not} linear.

\begin{prop}
The functors $\varphi_{X,V*}$ and $\varphi_{X,V}^\dagger$ induce an equivalence between the category of overconvergent crystals on $(X, V)$  with differential operators and the category of $i_X^{-1}\mathcal O_{V}$-modules and differential operators.
\end{prop}

\textbf {Proof : }
If $d : \mathcal F \to \mathcal G$ is a differential operator of order at most $n$, we can apply the functor $\varphi^\dagger_{X,V}$ to $\overline d$ and get a morphism of crystals
$$
\varphi^\dagger_{X,V}(\overline d) : \varphi^\dagger_{X,V}(p^{(n)}_{1*}p^{(n)\dagger}_{2} \mathcal F) \to \varphi^\dagger_{X,V}(\mathcal G).
$$
We have thanks to Proposition \ref{comdag} and \ref{finm},
$$
p^{(n)}_{1\AN^\dagger*}p^{(n)^{-1}}_{2\AN^\dagger}\varphi^\dagger_{X,V} \mathcal F = Êp^{(n)}_{1\AN^\dagger*}\varphi^\dagger_{X,V^{(n)}}(p^{(n)\dagger}_{2} \mathcal F) \simeq \varphi^\dagger_{X,V}(p^{(n)}_{1*}p^{(n)\dagger}_{2} \mathcal F).
$$
And we compose on the left with $p^{(n)^{*}}_{2\AN^\dagger}$ in order to get a differential operator at the crystal level.
Since we already know that our functor induce an equivalence between crystals on $(X, V)$  and of $i_X^{-1}\mathcal O_{V}$-modules when we consider only linear maps, the conclusion is immediate.
$\Box$

The following lemma shows that the category of overconvergent modules and differential operators has enough injectives.

\begin{lem}
Any complex $(E^\bullet, d)$ of overconvergent modules and differential operators has a resolution $(I^\bullet, d)$ made of injective overconvergent modules and differential operators.
\end{lem}

\textbf {Proof : }
It is sufficient to show that if $E \hookrightarrow I$ and $E' \hookrightarrow I'$ are two inclusion maps into injective overconvergent modules, then any is a differential operator $d : E \to E'$ of order at most $n$, extends to a differential operator $I \to I'$ of the same order.
By functoriality, it is clearly sufficient to extend
$$
\overline d : p^{(n)}_{1\AN^\dagger*}p^{(n)^{-1}}_{2\AN^\dagger} E \to E'.
$$
Since $I'$ is an injective module, this will follow from fact that both $p^{(n)}_{1\AN^\dagger*}$ and $p^{(n)^{-1}}_{2\AN^\dagger}$ are left exact.
$\Box$

Now that we can pull a complex of differential operators to $(X, V)_{\AN^\dagger}$, we want to push it to $(X/C, O)_{\AN^\dagger}$.
We first need an explicit description of this functor analog to Proposition \ref{restr}.

\begin{lem} \label{restr2}
Let 
$$
j_{X, V} : (X, V)_{\AN^\dagger} \to (X/C,O)_{\AN^\dagger}
$$
be the restriction morphism.
If $\mathcal F$ is an analytic sheaf on $(X, V)$ and $(X', V')$ is an analytic variety over $X/(C,O)$, we have
$$
(j_{X,V*} \mathcal F)_{X',V'} = p_{1*}\mathcal F_{X', V' \times_{O} V}.
$$
If $\mathcal F$ is a complex of abelian sheaves, we even get
$$
(Rj_{X,V*} \mathcal F)_{X',V'} = Rp_{1*}\mathcal F_{X', V' \times_{O} V}.
$$
\end{lem}

\textbf {Proof : }
The first assertion follows from the fact that $j_{X,V}^{-1}(X', V') = (X' , V' \times_{O} V)$ and therefore
$$
\Gamma((X', V'), j_{X,V*} \mathcal F) = \Gamma(j_{X,V}^{-1}(X', V'), \mathcal F) = \Gamma((X' , V' \times_{O} X), \mathcal F) $$
$$
= \Gamma(]X'[_{V' \times_{O} V}, \mathcal F_{X', V' \times V}) =  \Gamma(]X'[_{V'}, p_{1*}\mathcal F_{X', V' \times_{O} V}).
$$
For the second one, we recall that the realization functors are exact and preserve injectives.
$\Box$

In order to justify the next definition, we also prove the following.

\begin{lem} \label{linear}
If $d : E \to E'$ is a differential operator between overconvergent modules on $(X, V)$, the morphism
$$
j_{X,V*}(d) : j_{X,V*}E \to j_{X,V*}E'
$$
is $\mathcal O_{X/C,O}^\dagger$-linear.
\end{lem}

\textbf {Proof : }
This follows from the definition of a differential operator and the fact that the diagram
$$
\def\dar[#1]{\ar@<2pt>[#1]\ar@<-2pt>[#1]}
\xymatrix{(X,V^{(n)}) \dar[r]^{p^{(n)}_{1}}_{p^{(n)}_{2}} & (X, V) \ar[r]^{j_{X,V}} & X/(C, O)}
$$
is commutative.
$\Box$

\begin{dfn} If $\mathcal F$ is a sheaf of $i_{X}^{-1}\mathcal O_{V}$-modules, the \emph{linearization} of $\mathcal F$ is
$$
L(\mathcal F) = j_{X,V*}\varphi_{X,V}^{\dagger} \mathcal F.
$$
If $d : \mathcal F \to \mathcal G$ is a differential operator of finite order, its \emph{linearization} is
$$
L(d) :  j_{X,V*}\varphi_{X,V}^{\dagger}(d) : L(\mathcal F) \to L(\mathcal G).
$$

If $(\mathcal F^\bullet, d)$ is a complex of $i_{X}^{-1}\mathcal O_{V}$-modules and differential operators, the \emph{derived linearization} of $\mathcal F$ is
$$
RL(\mathcal F^\bullet) = Rj_{X,V*}\varphi_{X,V}^{\dagger} \mathcal F^\bullet.
$$

\end{dfn}

Note that $RL$ is not the right derived functor of $L$ (and anyway that $L$ is not left exact in general).
Thus, we have to be careful when one moves from $L$ to $RL$.
It also follows from lemma \ref{linear} that the linearization of a differential operator \emph{is} linear.
We now give a more explicit description of this linearization.

\begin{prop} \label{desc}
If $\mathcal F$ is a sheaf of $i_{X}^{-1}\mathcal O_{V}$-modules and $(X', V')$ is an analytic variety over $X/(C,O)$, we have
$$
L(\mathcal F)_{X', V'} = p_{1*}p_{2}^{\dagger}\mathcal F
$$
where $p_{1} : V' \times_{O} V \to V'$ and $p_{2} :  V' \times_{O} V \to V$ denote the projections.
Actually, if $(\mathcal F^\bullet, d)$ is a complex of $i_{X}^{-1}\mathcal O_{V}$-modules and differential operators, we have
$$
RL(\mathcal F^\bullet)_{X', V'} = Rp_{1*}p_{2}^{\dagger}\mathcal F^\bullet
$$
\end{prop}

\textbf {Proof : }
Results from lemma \ref{restr2}.
$\Box$

If $d : \mathcal F \to \mathcal G$ is a differential operator, we can also give an explicit description of the corresponding linear map
$$
p_{1*}p_{2}^{\dagger}\mathcal F \to p_{1*}p_{2}^{\dagger}\mathcal G
$$
on $]X'[_{V' \times_{O} V}$.
Of course, the main point is to describe the pull back morphism
$$
\delta : p_{1*}p_{2}^{\dagger}\mathcal F \to p_{1*}p_{2}^{\dagger}p_{2}^{(n)\dagger}\mathcal F
$$
that will be composed on the right with $\overline d$.
Thus, we consider the following diagram (where all products are taken relative to $O$)
$$
\def\dar[#1]{\ar@<2pt>[#1]\ar@<-2pt>[#1]}
\xymatrix{
V' \times V^{(n)} \ar@{^{(}->}[r]^{\Id\times i_{n}} \ar[d]_{p_{2}} & V' \times V \times V \dar[r]^{p_{13}}_{p_{12}} \ar[d]^{p_{2}} \dar[dr]^{p_{3}}_{p_{2}} & V' \times V \ar[r]^{p_{1}} \ar[d]^{p_{2}} & V'
\\ V^{(n)} \ar@{^{(}->}[r]^{i_{n}} & V \times V \dar[r]^{p_{2}}_{p_{1}} &  V & .}
$$
and call $p_{k}^{(n)} = p_{k} \circ i_{n}$ and $p_{1k}^{(n)} = p_{1k} \circ (\Id \times i_{n})$.
We consider the map
$$
p_{13}^{(n)*} : p_{2}^\dagger\mathcal F \to p_{13*}^{(n)}p_{13}^{(n)\dagger} p_{2}^\dagger\mathcal F = p_{13*}^{(n)} p_{2}^\dagger p_{2}^{(n)\dagger} \mathcal F
$$
and apply $p_{1*}$ in order to get
$$
\delta = p_{13}^{(n)*} : p_{1*}p_{2}^\dagger\mathcal F \to  p_{1*} p_{2}^\dagger p_{2}^{(n)\dagger} \mathcal F.
$$

\begin{prop} \label{poin} If $\mathcal F$ is a sheaf of $i_{X}^{-1}\mathcal O_{V}$-modules, there is a canonical isomorphism
$$
p_{X/C,O*}L(\mathcal F) \simeq u_{*} \mathcal F
$$
Actually, if $(\mathcal F^\bullet, d)$ is a complex of $i_{X}^{-1}\mathcal O_{V}$-modules and differential operators, we also have
$$
Rp_{X/C,O*}RL(\mathcal F^\bullet) \simeq Ru_{*} \mathcal F^\bullet.
$$
\end{prop}

\textbf {Proof : }
We consider the commutative diagram
$$
\xymatrix{
(X, V)_{\AN^\dagger} \ar[rr]^{\varphi_{X,V}} \ar[d]^ {j_{X,V}}
&& (]X[_{V})_{\an} \ar[d]^{u}
\\ (X/C,O)_{\AN^\dagger}  \ar[rr]^{p_{X/C,O}}
&& (]C[_{O})_{\an}
}
$$
We have
$$
p_{X/C,O*}L(\mathcal F) = p_{X/C,O*}j_{X,V*}\varphi_{X,V}^{\dagger} \mathcal F
=  u_{*} \varphi_{X,V*}\varphi_{X,V}^{\dagger} \mathcal F = u_{*} \mathcal F
$$
because $ \varphi_{X,V*}\varphi_{X,V}^{\dagger}\mathcal F = \mathcal F$.

Similarly, when $(\mathcal F^\bullet, d)$ is a complex of $i_{X}^{-1}\mathcal O_{V}$-modules and differential operators, we have
$$
Rp_{X/C,O*}RL(\mathcal F ^\bullet) = Rp_{X/C,O*}Rj_{X,V*}\varphi_{X,V}^{\dagger} \mathcal F ^\bullet
$$
$$
=  R(p_{X/C,O} \circ j_{X,V})_{*}\varphi_{X,V}^{\dagger} \mathcal F ^\bullet = R(u \circ  \varphi_{X,V})_{*}\varphi_{X,V}^{\dagger} \mathcal F ^\bullet
$$
$$
=  Ru_{*} \varphi_{X,V*}R\varphi_{X,V}^{\dagger} \mathcal F ^\bullet = Ru_{*} \varphi_{X,V*}\varphi_{X,V}^{\dagger} \mathcal F ^\bullet = Ru_{*} \mathcal F ^\bullet. \quad \Box
$$

%
%
\section{Grothendieck topology and overconvergence}

In Proposition \ref{locpoin} below, we will need to localize with respect to the Grothendieck topology.
Unfortunately, the dictionary of section 1.3 of \cite{Berkovich93} is not sufficient for us and we need to expand it a little bit.

We start with a geometrical result.

\begin{lem} \label{locpoint}
Let $(X \subset P)$ be formal embedding with $P$ affine and such that $\overline X\setminus X$ is a hypersurface $\overline g = 0$ in $\overline X$, with $g$ a function on $P$.
Let $V$ be an affinoid variety and $\lambda : V \to ]\overline X[_{P} \subset P_{K}$.
Then, the affinoid domains
$$
\{x \in V, \quad \epsilon |g(\lambda(x))| \geq 1 \}
$$
for $\epsilon \stackrel < \to 1$ form a cofinal family of neighborhoods of $X$ in $V$.
\end{lem}

\textbf {Proof : }
It is sufficient to remark that, by definition, we have
$$
]X[_{V} := \{x \in V, \overline g(\lambda(x)) \neq 0 \} = \{x \in V, |g(\lambda(x))| \geq 1 \}. \quad \Box
$$

If $V$ is an analytic variety over $K$, then $V_{G}$ will denotes the same set with its Grothendieck topology and structural sheaf.
If $u : V' \to V$ is a morphism of analytic varieties over $K$, then $u_{G} : V'_{G} \to V_{G}$ will denote the corresponding morphism of ringed sites.
Note that there is an obvious natural morphism of ringed sites $\pi_{V} : V_{G} \to V$.

If $(X, V)$ is an analytic variety with $]X[_{V}$ closed in $V$ and $\mathcal F$ (resp. $\mathcal G$) is a sheaf on $V$ (resp. $V_{G}$), one defines
$$
j^\dagger_{X} \mathcal F = \varinjlim j_{*}'j'^{-1} \mathcal F \quad \mathrm{and} \quad j_{X}^\dagger \mathcal G = \varinjlim j_{G*}'j_{G}'^{-1} \mathcal G
$$
where $j' : V' \hookrightarrow V$ runs through all embeddings of open neighborhoods of $X$ in $V$.
Recall from proposition 7.1 of \cite{LeStum06*} that $j_{X}^\dagger \mathcal F = i_{X*}i_{X}^{-1}\mathcal F$.

\begin{dfn} An analytic variety $(X, V)$ over $\mathcal V$ is said
\begin{enumerate}
\item\emph{paracompact} (resp. \emph{countable at infinity}), if there exists an open neighborhood $V'$ of $X$ in $V$ which is paracompact (resp. countable at infinity) with $]X[_{V}$ closed in $V'$.
\item \emph{Hausdorff} (resp. \emph{separated}) if there exists a neighborhood $V'$ of $X$ in $V$ which is Hausdorff (resp. separated).
\end{enumerate}
\end{dfn}

\begin{lem} \label{comolim}
If $(X, V) $ is a separated analytic variety over $\mathcal V$ that is countable at infinity and $\mathcal G$ is a coherent sheaf on $V_{G}$, then
$$
\forall q \geq 0, \quad \varinjlim H^q(V'_{G}, \mathcal G) \simeq H^q(V_{G}, j_{X}^\dagger \mathcal G)
$$
where $V'$ runs through all the open neighborhoods of $X$ in $V$.
\end{lem}

The result actually holds for any (good) paracompact analytic variety (this follows from Proposition \ref{compov} below).

\textbf {Proof : }
Note first that this equality is always true when $V$ is compact because then, $V_{G}$ is quasi-compact and quasi-separated. In particular, this holds in the affinoid case. If, moreover, $X$ has a cofinal family of affinoid neighborhoods in $V$, we see that the left hand side is zero for $q > 0$ and it follows that the right-hand side too is zero.

After blowing up the locus at infinity of $X$ and localizing on $P$, we see that any analytic variety has an affinoid covering whose terms are as in lemma \ref{locpoint}.
It follows that, in our situation, there exists a countable admissible affinoid covering of finite type $V = \cup_{i\in \mathbf N} W_{i}$ and for each $i \in \mathbf N$ a decreasing family $\{W_{i,\epsilon}\}_{\epsilon \stackrel < \to 1}$ of affinoid neighborhoods of $X$ in $W$ that are cofinal.
If $\underline \epsilon := \{\epsilon_{i}\}_{i \in \mathbf N}$ is any sequence inside $[0, 1[$, we let $V_{\underline \epsilon} := \cup W_{i, \epsilon_{i}}$.
This is a  covering of finite type by affinoid domains and therefore an admissible covering of an analytic domain.
Note that $V_{\underline \epsilon}$ is a neighborhood of $X$ in $V$ because this question is local for the Grothendieck topology of $V$ and that $V_{\underline \epsilon} \cap W_{i} \supset W_{i,\epsilon_{i}}$
Moreover, if $V'$ is an open neighborhood of $X$ in $V$, then, for each $i \in \mathbf N$, $V' \cap W_{i}$ is an open neighborhood of $X$ in $W_{i}$ and there exits $\epsilon_{i}$ such that $W_{i,\epsilon_{i}} \subset V'$. It follows that $V_{\underline \epsilon} \subset V'$.

For fixed $i_{1}, \ldots, i_{p} \in \mathbf N$, the affinoid subsets $W_{\underline {i, \epsilon}} := W_{i_{1},\epsilon_{1}} \cap \cdots \cap W_{i_{p},\epsilon_{p}}$ with $\epsilon_{1}, \ldots, \epsilon_{p}  < 1$, form a cofinal system of neighborhoods of $X$ in $W_{\underline i} := W_{i_{1}} \cap \cdots \cap W_{i_{p}}$.
Since our variety is separated, all the analytic domains $W_{\underline {i, \epsilon}}$ are all affinoid, and we have for all $q > 0$,
$$
H^q(W_{\underline {i, \epsilon},G}, \mathcal G) = 0 \quad \mathrm{and} \quad H^q(W_{\underline i,G}, j_{X}^\dagger\mathcal G) = 0
$$
(use the remark at the beginning of this proof).
It follows that for any sequence $\underline \epsilon$, we have
$$
H^q(V_{\underline \epsilon ,G}, \mathcal G) = \check H^q(\{W_{i,\epsilon_{i},G}\}, \mathcal G) \quad \mathrm{and} \quad H^q(V,j_{X}^\dagger \mathcal G) = \check H^q(\{W_{i,G}\},j_{X}^\dagger \mathcal G))
$$
Since filtered direct limits are exact, it is therefore sufficient to show that
$$
\forall q \geq 0, \quad \varinjlim_{\underline \epsilon} \check C^q(\{W_{i,\epsilon_{i,G}}\}, \mathcal G) \simeq \check C^q(\{W_{i,G}\}, j_{X}^\dagger \mathcal G).
$$
We are therefore reduced to prove that for all $i_{1}, \ldots, i_{p} \in \mathbf N$, we have
$$
\varinjlim_{\underline \epsilon} \prod_{\underline i}\Gamma(W_{\underline{i,\epsilon},G}, \mathcal G) \simeq \prod_{i}\varinjlim_{\epsilon}\Gamma(W_{\underline i,\epsilon, G}, \mathcal G)
$$
which is an easy exercise on commutation of products and filtered direct limits.
$ \Box$

\begin{prop} \label{compov}
If $(X, V)$ is a (good) analytic variety on $\mathcal V$ with $]X[_{V}$ closed in $V$ and $\mathcal F$ is an $\mathcal O_{V}$-module, we have $\pi_{V*}j_{X}^\dagger\pi_{V}^*\mathcal F = j^\dagger \mathcal F$.
Actually, if $\mathcal F$ is coherent, then
$$
R\pi_{V*}j_{X}^\dagger\pi_{V}^*\mathcal F = j^\dagger \mathcal F.
$$
\end{prop}

\textbf {Proof : }
The first assertion is checked on the stalks where it is clear since any point has a basis of affinoid neighborhoods.
In order to prove the second one, since any good analytic variety has a basis of open subsets that are separated and countable at infinity, it is sufficient to show that, if $V$ is separated and countable at infinity, then
$$
\forall q \geq 0, \quad H^q(V_{G}, j_{X}^\dagger\pi_{V}^*\mathcal F) = H^q(V, j_{X}^\dagger \mathcal F).
$$
We proved in lemma \ref{comolim} that
$$
\forall q \geq 0, \quad \varinjlim H^q(V'_{G}, \pi_{V}^*\mathcal F) \simeq H^q(V_{G}, j_{X}^\dagger \pi_{V}^*\mathcal F)
$$
where where $V'$ runs through all the open neighborhoods of $X$ in $V$.
On the other side, since $V$ is paracompact and $]X[_{V}$ is closed in $V$, we also have (see for example Remark 2.6.9 of \cite{KashiwaraSchapira90})
$$
\forall q \geq 0, \quad \varinjlim H^q(V', \mathcal F) \simeq H^q(V, j_{X}^\dagger \mathcal F).
$$
Our assertion therefore follows from Proposition 1.3.6.ii of \cite{Berkovich93}.
$\Box$

In fact, we will need a relative version of this result:

\begin{cor} \label{Grot}
Let $(X', V') \to (X, V)$ be a morphism of (good) analytic varieties with $]X[_{V}$ closed in $V$ (resp. $]X[_{V'}$ closed in $V'$).
Let $\mathcal F$ (resp. $\mathcal F'$) be a complex with coherent terms on $V$ (resp. $V'$) and differentials defined on $V_{G}$ (resp. $V_{G'}$).
Assume that there is an isomorphism $j^\dagger_{X}\pi_{V}^*\mathcal F \simeq Ru_{G*} j_{X'}^\dagger \pi_{V'}^*\mathcal F'$.
Then, we also have $j^\dagger_{X}\mathcal F \simeq Ru_{*} j^\dagger_{X'}\mathcal F'$.
\end{cor}

\textbf {Proof : }
We apply $R\pi_{V*}$ on both sides of the first equality in order to get, thanks to lemma \ref{compov},
$$
j^\dagger_{X}\mathcal F \simeq R\pi_{V*}j^\dagger_{X}\pi_{V}^*\mathcal F \simeq R\pi_{V*}Ru_{G*} j_{X'}^\dagger \pi_{V'}^*\mathcal F'
$$
$$
\simeq Ru_{*}R\pi_{V'*} j_{X'}^\dagger \pi_{V'}^*\mathcal F' = Ru_{*} j_{X'}^\dagger \mathcal F' . \quad \Box
$$

%
%
\section{Cohomology}

We assume in this last section that $\Char K = 0$.

Let $(X, V) \to (C, O)$ be a morphism of analytic varieties over $\mathcal V$.
If $\mathcal F$ is an $i_{X}^{-1}\mathcal O_{V}$-module with an integrable connexion, one can build its de Rham complex
$$
\mathcal F \otimes_{i_{X}^{-1}\mathcal O_{V} } i_{X}^{-1}\Omega^\bullet_{V/O} 
$$
(we will only consider the smooth case).
It is a complex of differential operators of order at most $1$.
We may therefore consider its derived linearization 
$$
RL(\mathcal F \otimes_{i_{X}^{-1}\mathcal O_{V} } i_{X}^{-1}\Omega^\bullet_{V/O}).
$$

We can give an explicit description of this object (we stick to the case of crystals):

\begin{lem} \label{resol}
Let $(X, V) \to (C, O)$ be a morphism of analytic varieties over $\mathcal V$ and $E$ an overconvergent crystal on $X/(C,O)$.
Then
\begin{enumerate}
\item If $(X', V')$ is an analytic variety over $X/(C,O)$, we have
$$
RL(E_{X,V} \otimes_{i_{X}^{-1}\mathcal O_{V} } i_{X}^{-1}\Omega^\bullet_{V/O})_{X',V'} = Rp_{1*}(p_{1}^\dagger E_{X',V'} \otimes_{i_{X'}^{-1}\mathcal O_{V' \times V} } i_{X'}^{-1}\Omega^\bullet_{V' \times V/V'}).
$$
\item
There is a canonical augmentation
$$
E \to RL(E_{X,V} \otimes_{i_{X}^{-1}\mathcal O_{V} } i_{X}^{-1}\Omega^\bullet_{V/O}).
$$
\end{enumerate}
\end{lem}

Actually, both results hold before deriving.

\textbf {Proof : }
Since $E$ is a crystal, we have
$$
p_{2}^{\dagger} (E_{X,V} \otimes_{i_{X}^{-1}\mathcal O_{V} } i_{X}^{-1}\Omega^\bullet_{V/O}) = p_{2}^{\dagger}E_{X,V} \otimes_{i_{X'}^{-1}\mathcal O_{V' \times V} } i_{X'}^{-1}\Omega^\bullet_{V' \times V/V'}
$$
$$
= E_{X',V' \times V} \otimes_{i_{X'}^{-1}\mathcal O_{V' \times V} } i_{X'}^{-1}\Omega^\bullet_{V' \times V/V'}
= p_{1}^\dagger E_{X',V'} \otimes_{i_{X'}^{-1}\mathcal O_{V' \times V} } i_{X'}^{-1}\Omega^\bullet_{V' \times V/V'}
$$
and we may apply $Rp_{1*}$ in order to get the first assertion thanks to Proposition \ref{desc}.

Note that it is sufficient to prove the second assertion before deriving and then compose with the canonical morphism $L \to RL$.
Since $E$ is a crystal, we have
$$
j^{-1}_{X,V}E = \varphi_{X,V}^{\dagger}\varphi_{X,V*}j^{-1}_{X,V}E = \varphi_{X,V}^{\dagger}E_{X,V}
$$
and by adjunction, we obtain a map
$$
E \to j_{X,V*}\varphi_{X,V}^{\dagger}E_{X,V} = L(E_{X,V}).
$$
We have to show that the composite map
$$
E \to  L(E_{X,V}) \to L(E_{X,V} \otimes_{i_{X}^{-1}\mathcal O_{V} } i_{X}^{-1}\Omega^1_{V/O})
$$
is zero.
This can be done locally and, thanks to the first part, we are reduced to check that
$$
E_{X', V'} \to p_{1*}p_{1}^\dagger E_{X',V'} \to p_{1*}(p_{1}^\dagger E_{X',V'} \otimes_{i_{X'}^{-1}\mathcal O_{V' \times V} } i_{X'}^{-1}\Omega^\bullet_{V' \times V/V'})
$$
is zero.
This is linear in $E_{X', V'}$ and we are reduced to show that the composite map
$$
i_{X'}^{-1}\mathcal O_{V'} \to p_{1*}i_{X'}^{-1}\mathcal O_{V' \times V} \to p_{1*} i_{X'}^{-1}\Omega^\bullet_{V' \times V/V'}
$$
is zero.
Here again, we need to be careful since the base change map $i_{X'}^{-1}p_{1*}\ \to p_{1*}i_{X'}^{-1}$ is not bijective in general.
But we do not care because it is \emph{sufficient} that the composition
$$
i_{X'}^{-1}\mathcal O_{V'} \to i_{X'}^{-1}p_{1*}\mathcal O_{V' \times V} \to i_{X'}^{-1}p_{1*}\Omega^\bullet_{V' \times V/V'}
$$
is zero.
Finally, we may drop the $i_{X'}^{-1}$ and are reduced to the augmentation from the constant sheaf to the usual de Rham complex.
$\Box$

\begin{thm} \label{derham}
Let
$$
\xymatrix{X \ar@{^{(}->}[r] \ar[d] & P \ar[d]^u & V \ar[l] \ar[d] \\ C \ar@{^{(}->}[r] & S & O \ar[l]}
$$
be a morphism of analytic varieties with $u$ proper and smooth at $X$, $P$ admissible and the left square cartesian.
Let $E$ be an overconvergent module of finite presentation on $X/(C, O)$.
Then, there is a canonical isomorphism
$$
Rp_{X/C,O*}E \simeq Ru_{*} (E_{X,V} \otimes_{i_{X}^{-1}\mathcal O_{V} } i_{X}^{-1}\Omega^\bullet_{V/O}).
$$
\end{thm}

\textbf {Proof : }
It is sufficient, thanks to Proposition \ref{poin}, to prove the Proposition \ref{poinlem} below.
$\Box$

\begin{prop} \label{poinlem}
In the situation of the theorem, the augmentation map is an isomorphism
$$
E \simeq RL(E_{X,V} \otimes_{i_{X}^{-1}\mathcal O_{V} } i_{X}^{-1}\Omega^\bullet_{V/O}).
$$
\end{prop}

\textbf {Proof : }
We have to show that for any analytic variety $(X', V')$ on $X/(C, O)$, we have
$$
E_{X', V'} \simeq RL(E_{X,V} \otimes_{i_{X}^{-1}\mathcal O_{V} } i_{X}^{-1}\Omega^\bullet_{V/O})_{X', V'}.
$$
If we write $\mathcal F' := E_{X', V'}$, lemma \ref{resol} allows us to rewrite this isomorphism as
$$
\mathcal F' \simeq Rp_{1*}(p_{1}^\dagger\mathcal F' \otimes_{i_{X'}^{-1}\mathcal O_{V' \times V} } i_{X'}^{-1}\Omega^\bullet_{V' \times V/V'}).
$$
Note that $p_{1} : P' \times P \to P'$ is a morphism of good embeddings of $X'$ that is proper and smooth at $X'$.
Moreover, we may replace $P'$ by its maximum flat formal subscheme $P'^{\mathrm{fl}}$ and consequently $X'$ with its trace on $P'^{\mathrm{fl}}$.
Our assertion therefore follows from Proposition \ref{locpoin} below.
$\Box$

\begin{prop} \label{locpoin}
Let $u : P' \to P$ be morphism of good admissible embeddings of $X$ that is proper and smooth at $X$.
Let $V$ be a good analytic variety, $\lambda : V \to P_{K}$ a morphism and $V' := P'_{K} \times_{P_{K}} V$.
If $\mathcal F$ is a coherent $i_{X}^{-1}\mathcal O_{V}$-module, then
$$
\mathcal F \simeq Ru_{*}(u^\dagger\mathcal F \otimes_{i_{X}^{-1}\mathcal O_{V'} } i_{X}^{-1}\Omega^\bullet_{V'/V}).
$$
\end{prop}

\textbf {Proof : }
It follows form corollary \ref{locx} that the question is local on $X$.
We want to show that it is local for the Grothendieck topology of $V$ also.
This is not clear.
First of all, since the question is local for the analytic topology, we may assume that $V$ is paracompact and replace, as we usually do, $\mathcal F$ by $i_{X}^{-1}\mathcal F$ where $\mathcal F$ is a coherent $\mathcal O_{V}$-module.
After pushing by $i_{X*}$, we are reduced, thanks to proposition 7.1 of \cite{LeStum06*}, to show that
$$
j_{X}^\dagger \mathcal F \simeq Ru_{*}(j_{X}^\dagger u^*\mathcal F \otimes_{\mathcal O_{V'} } \Omega^\bullet_{V'/V}).
$$
Since $\mathcal F$ is coherent and $V$ and $V'$ are good (use Proposition 2.3 of \cite{LeStum06*} for $V'$), it is sufficient to prove thanks to lemma \ref{Grot} that
$$
j_{X}^\dagger \pi_{V}^*\mathcal F \simeq Ru_{G*}(j_{X}^\dagger u_{G}^*\pi_{V}^*\mathcal F \otimes_{\mathcal O_{V_{G}'} } \Omega^\bullet_{V'/V,G}).
$$

Now, the question is local for the Grothendieck topology of $V$.
Thanks to the strong fibration theorem \ref{strong}, we may therefore assume that $P$ is affine, that the locus at infinity of $X$ is a hypersurface, that $V$ is affinoid and contained inside $]\overline X[_{V}$, and that $P' = \widehat {\mathbf A}^{n}_{P}$.

One may then proceed by induction on $n$.
The case $n = 1$ results from lemma \ref{smallloc} below.
For the induction process, we use the Gauss-Manin connexion as in lemma 6.5.5 of \cite{LeStum07*}.
$\Box$

\begin{lem} \label{smallloc}
Let $(X \subset P)$ be a formal embedding with $P$ affine and such that the locus at infinity of $X$ is a hypersurface.
Let $\lambda : V \to ]\overline X[_{P} \subset P_{K}$ be a morphism with $V$ affinoid.
If $\mathcal F$ is a coherent $i_{X}^{-1}\mathcal O_{V}$-module, there is canonical isomorphism
$$
\Gamma(]X[_{V}, \mathcal F) \simeq R\Gamma(]X[_{V} \times \mathbf D(0, 1^{-}),  p_{1}^\dagger \mathcal F \stackrel {\partial/\partial t} \longrightarrow  p_{1}^\dagger \mathcal F).$$
\end{lem}

\textbf {Proof : }
This is proved exactly as in proposition 6.5.7 of \cite{LeStum07*}.
$\Box$

We may now derive several corollaries.

\begin{cor}
Let $(C, O)$ be an analytic variety over $\mathcal V$, $f : X' \to X$ a morphism of algebraic varieties over $C$, $E'$ an overconvergent module of finite presentation on $X'/(C,O)$ and $(U \subset P \leftarrow V)$ be an analytic variety over $X/(C, O)$.

Let $U' := U \times_{X} X'$ and
$$
\xymatrix{U' \ar@{^{(}->}[r] \ar[d] & P' \ar[d]^u & V' \ar[l] \ar[d] \\ U \ar@{^{(}->}[r] & P & V \ar[l]}
$$
be a commutative diagram with $u$ proper and smooth at $U'$, $P'$ admissible and the left square cartesian.
Then, there is a canonical isomorphism
$$
(Rf_{\AN^\dagger*}E')_{(U, V)} \simeq Ru_{*} (E'_{U',V'} \otimes_{i_{U'}^{-1}\mathcal O_{V'} } i_{U'}^{-1}\Omega^\bullet_{V'/V})
$$
\end{cor}

\textbf {Proof : }
We have
$$
(Rf_{\AN^\dagger*}E)_{(U, V)} = Rp_{U \times_{X} X'/(U, V)*} E_{|U \times_{X} X'/(U, V)}
$$
and we apply theorem \ref{derham}.
$\Box$

If $V$ is a Hausdorff analytic variety over $K$, then the set $V_{rig}$, or $V_{0}$ for short, of rigid points of $V$ inherits a structure of rigid analytic variety from the Grothendieck topology of $V$ and there is a morphism of sites $\pi_{0} : V_{0} \to V$.
Actually, we get an equivalence of toposes $\widetilde{V_{0}} \simeq \widetilde{V_{G}}$.
If $\mathcal F$ is a complex of abelian groups on $V_{0}$, we will write $\mathcal F^{\an} := R\pi_{0*}\mathcal F$.
\begin{prop} \label{rigcomp}
Let $S$ be a formal $\mathcal V$-scheme and $p: X \to C$ a morphism of algebraic varieties over $S_{k}$.
Let $E$ be an overconvergent isocrystal on $X/S$.
Let $C \hookrightarrow Q$ be a good formal embedding over $S$.
Then,  we have
$$
(Rp_{\rig} E)^{\an}= i_{C*}Rp_{X/C\subset Q*} E
$$
\end{prop}

Of course, we use Proposition 9.4 of \cite{LeStum06*} to identify the category of overconvergent modules of finite presentation and the category of overconvergent isocrystals.

\textbf {Proof : }
In order to define rigid cohomology, we need to assume that $p$ extends to a morphism $u : P \to Q$ that is proper and smooth at $X$ (in general, we would need some glueing which is possible as shown in the appendix \ref{Zaris}).
Now, we see $E$ as an overconvergent module and consider its realization $E_{P}$ on $P$.
The question being local on $]C[_{Q}$, we may assume that $Q$ is quasi-compact and therefore also that $P$ is quasi-compact.
Then, there exists a good open neighborhood $V$ of $]X[_{P}$ in $]\overline X[_{P}$ and a coherent module with an integrable connexion $\mathcal F$ on $V$ such that $E_{P} = i_{X}^{-1} \mathcal F$.
Since $V$ is good, $\mathcal F$ extends uniquely to a coherent module with an integrable connexion $\mathcal F_{0}$ on $V_{0}$.
Recall from proposition 4.1 of \cite{LeStum06*} that strict neighborhoods in rigid geometry correspond bijectively to neighborhoods in Berkovich theory.
In particular, if we let $W_{0} := ]\overline C[_{Q}$, we have, by definition,
$$
Rp_{\rig} E := Ru_{0*} (j_{0}^\dagger \mathcal F_{0} \otimes_{\mathcal O_{V_{0}}} \Omega^\bullet_{V_{0}/W_{0}})
$$
Since $V$ is good, it follows from Lemma \ref{compov} (and the equivalence $\widetilde{V_{0}} \simeq \widetilde{V_{G}}$) that
$$
(j_{0}^\dagger \mathcal F_{0} \otimes_{\mathcal O_{V_{0}}} \Omega^\bullet_{V_{0}/W_{0}})^{\an} = j^\dagger \mathcal F \otimes_{\mathcal O_{V}} \Omega^\bullet_{V/W} = i_{X*}i_{X}^{-1}\mathcal F \otimes_{\mathcal O_{V}} \Omega^\bullet_{V/W}
$$
and therefore,
$$
(Rp_{\rig} E)^{\an} = Ru_{*} (i_{X*}i_{X}^{-1}\mathcal F \otimes_{\mathcal O_{V}} \Omega^\bullet_{V/W}) = Ru_{*} i_{X*}(i_{X}^{-1}\mathcal F \otimes_{i_{X}^{-1}\mathcal O_{V}} i_{X}^{-1}\Omega^\bullet_{V/W}) 
$$
$$
= i_{C*}Ru_{*} (E_{P} \otimes_{i_{X}^{-1}\mathcal O_{V}} i_{X}^{-1}\Omega^\bullet_{V/W}) = i_{C*}Rp_{X/C\subset Q*} E. \quad \Box
$$

As a corollary, we recover the fact that the rigid cohomology of an overconvergent isocrystal on $X/S$ is independent of the choice of the formal embedding $X \hookrightarrow P$ into a formal scheme that is proper and smooth at $X$ over $S$.

There is a particular case of the proposition that is worth stating:

\begin{cor}
If $X$ is an algebraic variety over $k$ and $E$ an overconvergent isocrystal on $X/K$, we have for all $i \in \mathbf N$,
$$
H^{i}_{\rig}(X/K, E) = H^{i}((X/\mathcal V)_{\AN^\dagger}, E).
$$
\end{cor}

$\Box$

%
%
\appendix
\section{Zariski localization}\label{Zaris}

In this section, we show that the cohomology of crystals on an analytic variety $(X, V)$ over $\mathcal V$ is local for the Zariski topology of $X$.
This is necessary if one wants to show that rigid cohomology coincide with our cohomology even when there is no global embedding that is proper and smooth at $X$.

We let $(C, O)$ be an analytic variety over $\mathcal V$ but we first give a general definition:

\begin{dfn}
Let $T$ be an analytic presheaf on $\mathcal V$.
An analytic sheaf $\mathcal F$ on $T$ is said \emph{of Zariski type} if for any analytic variety $(X,V)$ over $T$ and any open immersion $\alpha : U \hookrightarrow X$, we have $]\alpha[_{V}^{-1}\mathcal F_{X,V} = \mathcal F_{U,V}$.
\end{dfn}

Clearly, this property will then be satisfied for any locally closed immersion.

Recall that if $X$ is an algebraic variety over $C$ and $\alpha : U \hookrightarrow X$ is a locally closed immersion, then
$$
\alpha_{\AN^\dagger} : (U/C,O)_{\AN^\dagger} \hookrightarrow (X/C,O)_{\AN^\dagger}
$$
denotes the corresponding morphism of toposes (this is a restriction map).

\begin{lem} Let  $X$ be an algebraic variety over $C$.
Let $\mathcal F$ be a sheaf of Zariski type on $X/(C,O)$ and $\alpha : U \hookrightarrow X$ be a locally closed immersion over $C$.
Let $(X', V')$ be an analytic variety over $X/(C,O)$, $f : X' \to X$ the structural map, $U' = f^{-1}(U)$ and $\alpha' : U' \hookrightarrow X'$ the inclusion map.
Then, we have
$$
(\alpha_{\AN^\dagger*}\alpha_{\AN^\dagger}^{-1}\mathcal F)_{X',V'} = ]\alpha'[_{*}]\alpha'[^{-1} \mathcal F_{X',V'}.
$$
\end{lem}

\textbf {Proof : }It follows from Proposition 5.9 of \cite{LeStum06*} (or Proposition \ref{restr} of this article) that
$$
(\alpha_{\AN^\dagger*}\alpha_{\AN^\dagger}^{-1}\mathcal F)_{X',V'} =  ]\alpha'[_{*}(\alpha_{\AN^\dagger}^{-1}\mathcal F)_{U',V'}
$$
$$
= ]\alpha'[_{*}\mathcal F_{U',V'} = ]\alpha'[_{*}]\alpha'[^{-1} \mathcal F_{X',V'}. \quad \Box
$$

\begin{prop}
Let  $X$ be an algebraic variety over $C$ and $\mathcal F$ an abelian sheaf of Zariski type on $X/(C,O)$.
If $\{X_{k}\}_{k \in I}$ is a locally finite open covering of $X$, there is a long exact sequence
$$
0 \to \mathcal F \to \prod_{k} \alpha_{k\AN^\dagger*}\alpha^{-1}_{k\AN^\dagger}  \mathcal F \to  \prod_{k} \alpha_{kl\AN^\dagger*}\alpha^{-1}_{kl\AN^\dagger}  \mathcal F \to \cdots 
$$
\end{prop}

\textbf {Proof : }The assertion can be checked on some $(X', V')$ with $f : X' \to X$. We may clearly assume that $X' = X$ and we are reduced, thanks to the previous lemma, to the exact sequence of Proposition \ref {Zar}.
$\Box$

Thus, we get a spectral sequence of absolute cohomology:

\begin{cor}
With the assumptions of the proposition, there is a spectral sequence
$$
E_{1}^{r,s} = \bigoplus_{|\underline k| = r + 1} R^sp_{X_{\underline k}/C,O*} \mathcal F_{|X_{\underline k}} \Rightarrow R^{r+s}p_{X_{V}/C,O*}\mathcal F
$$
where
$$
X_{\underline k} := X_{k_1} \cap \cdots X_{k_{r+1}}.
$$
\end{cor} $\Box$

And we finally mention the spectral sequence of relative cohomology:

\begin{cor}
Let  $f : X' \to X$ be a morphism of analytic varieties over $C$, $\mathcal F'$ an abelian sheaf of Zariski type on $X'/(C,O)$, $\{X'_{k}\}_{k \in I}$ a locally finite open covering of $X'$ and $f : X' \to X$ any $C$-morphism.
Then, there is a spectral sequence
$$
E_{1}^{r,s} = \bigoplus_{|\underline k| = r + 1} R^sf_{|X'_{\underline k}\AN^\dagger*}\mathcal F'_{|X'_{\underline k}} \Rightarrow R^{r+s}f_{\AN^\dagger*}\mathcal F'.
$$
\end{cor} $\Box$

The important thing is that the above results apply to crystals as the following proposition shows.
Moreover, we also see that, being a crystal is of local nature for the Zariski topology.

\begin{prop} \label{crisloc} Let $X$ be an algebraic variety over $C$ and $E$ an overconvergent module on $X/(C,O)$.
Let $\{X_{k}\}_{k \in I}$ be a locally finite open covering of $X$.
Then $E$ is a crystal if and only if it is of Zariski type and for each $k$, $E_{|X_{k}}$ is a crystal.
\end{prop}

\textbf {Proof : }
It follows from the second assertion of Proposition 5.9 of \cite{LeStum06*} that an analytic crystal is of Zariski type.
And we also know that the inverse image of a crystal is a crystal.

Conversely, assume that $E$ is an overconvergent module of Zariski type and that for each $k$, $E_{|X_k}$ is a crystal.
Let $(f, u) : (U', V') \to (U, V)$ be any morphism of analytic varieties over $X$ and $g : U \to X$ the canonical map.
For each $k = (k_1, \ldots, k_r)$, let
$$
X_{\underline k} = X_{k_1} \cap \cdots \cap X_{k_r}, \quad U_{\underline k} = g^{-1}(X_{\underline k}), \quad U'_{\underline k} = f^{-1}(U_{\underline k}),
$$
$f_{\underline k} : U'_{\underline k} \to U_{\underline k}$ the restriction of $f$ and finally, $\alpha_k : U_{\underline k} \hookrightarrow U, \alpha'_{\underline k} : U'_{\underline k} \hookrightarrow U'$ the inclusion maps.
Since $E_{|X_{\underline k}}$ is a crystal, we have $]f_{\underline k}[_u^\dagger E_{U_{\underline k},V} = E_{U'_{\underline k},V'}$. Since $E$ is Zariski type, it follows that
$$
]\alpha'_{\underline k}[^{-1} ]f[_u^\dagger E_{U,V} = ]f_{\underline k}[_u^\dagger]\alpha_{\underline k}[^{-1} E_{U,V} = ]\alpha'_{\underline k}[^{-1} E_{U',V'}
$$
and we can apply Lemma \ref{Zar}.
$\Box$

%
%

\section{\v{C}ech-Alexander cohomology}

It is sometimes necessary to compute the cohomology of general abelian sheaves that are not crystals. This would be important for example if one is willing to develop a Dieudonn\'e crystalline theory on the overconvergent site.

We will let as usual $(X, V) \to (C, O)$ be a morphism of analytic varieties over $\mathcal V$ but we start with general considerations.
Let $C$ be a standard site and $T$ an object of $C$. From the simplicial complex $(T^{r+1})_{r \in \mathbf N}$, we can form for any abelian sheaf $\mathcal F$ on $C$, a complex $\mathcal C_T(\mathcal F)$ whose terms are given by
$$
\mathcal C^r_T(\mathcal F) = j_{r*}j_{r}^{-1}\mathcal F
$$
where $j_{r} :  \mathcal C_{/T^{r+1}} \to C$ is the canonical map. There is a canonical augmentation $\mathcal F \to \mathcal C_T(\mathcal F)$. This is a quasi-isomorphism if $T$ is a covering of the final object of $\tilde C$. 
Actually, since this construction is functorial, we get for any complex of abelian sheaves $\mathcal F$, a quasi-isomorphism $\mathcal F \simeq [C(\mathcal F)]$ where $[E]$ denotes the complex associated to the a bicomplex $E$.
Finally, since both $j_{r*}$ and $j_{r}^{-1}$ both preserve injective, we see that if $\mathcal I$ is a  complex of injective sheaves, so is $[C(\mathcal I)]$.

Now, we let $X_V/C,O)$ be the image of $(X, V)$ inside the presheaf $X/(C, O)$ and $\AN^\dagger(X_V/C,O)$ the corresponding site (objects are in $\AN^\dagger(X, V)$ and morphisms in $\AN^\dagger(X/C,O)$).
For each $r \in \mathbf N$, we can embed $X$ diagonally into
$$
V^{r+1} := V \times_{O} \times \cdots \times_{O} V.
$$
Let $\mathcal F$ be an analytic sheaf on $\AN^\dagger(X_V/C,O)$. From the projections and degeneracy maps
$$
d^r_i : V^{r+1} \to V^r \textrm{ and } s^r_i : V^{r-1} \to V^r,
$$
we derive morphisms
$$
(d^r_i)^{-1}\mathcal F_{V^{r}} \to \mathcal F_{V^{r+1}} \textrm{ and } (s^r_i)^{-1}\mathcal F_{V^{r}} \to \mathcal F_{V^{r-1}},
$$
or by adjunction
$$
\mathcal F_{V^{r}} \to d^r_{i*} \mathcal F_{V^{r+1}} \textrm{ and } \mathcal F_{V^{r}} \to d^r_{i*}\mathcal F_{V^{r-1}}.
$$
Now if $p_{1} : V^{r+1} \to V$ denotes the first projection, we may apply $p_{1*}$ and get a cosimplicial complex $( p_{1*}\mathcal F_{V^{r}}, \delta_i^r, \sigma_i^r)$.

\begin{dfn}
With the above notations, the \emph{Cech-Alexander complex} $\mathcal {CA}(\mathcal F)$ of $\mathcal F$ is the complex on $]X[_V$ which is formally derived from this cosimplicial complex.

When $\mathcal F$ is a complex of abelian sheaves, the  \emph{derived Cech-Alexander complex} of $\mathcal F$ is the complex $R\mathcal {CA}(\mathcal F)$ associated to the bicomplex $\mathcal {CA}(\mathcal I)$ where $\mathcal I$ is a resolution of $\mathcal F$.
\end{dfn}

Concerning $\mathcal {CA}(\mathcal F)$, we see that
$$
d^r = \sum_{i = 0}^r (-1)^i \delta_i^r : p_{1*}\mathcal F_{V^{r}} \to  p_{1*}\mathcal F_{V^{r+1}}.
$$
The degeneracy arrows are used to check that this is a complex.
All this is clearly functorial.

We denote by
$$
p_{V^{r+1}} : ]X[_{V^{r+1}} \to ]C[_{O}
$$
the structural map and by
$$
p_{X_V/C,O} : (X_V/C,O)_{\AN^\dagger} \to (]C[_O)_{\an}
$$
the canonical map.

\begin{prop}
Let $(X, V) \to (C, O)$ be a morphism of analytic varieties over $\mathcal V$.
If $\mathcal F$ is an abelian sheaf on $X_V/C,O$, we have
$$
p_{X_V/C,O*}\mathcal F = p_{V*} \mathcal {CA}(\mathcal F).
$$
If $\mathcal F$ is a complex of abelian sheaves on on $X_V/C,O$, we have an isomorphism
$$
Rp_{X_V/C,O*}\mathcal F = Rp_{V*} R\mathcal {CA}(\mathcal F).
$$
\end{prop}

\textbf {Proof : }
We can apply the first considerations at the beginning of this section to our situation. 
Namely, by definition, $(X, V)$ is a covering of the final object of $(X_V/C,O)_{\AN^\dagger}$.
Thus, any abelian sheaf $\mathcal F$ on $X_V/C,O$ has a canonical resolution $\mathcal F \to \mathcal C(\mathcal F)$ with $\mathcal C^r(\mathcal F) =  j_{r*}j_{r}^{-1}\mathcal F$ where
$$
j_{r} : (X, V^{r+1})_{ \AN^\dagger} \to X_V/C,O)_{ \AN^\dagger}
$$
is the restriction map.

For each $r$, we have a commutative diagram
$$
\xymatrix{
(X, V^{r+1})_{\AN^\dagger} \ar[rr]^{\varphi_{V^{r+1}}} \ar[d]^{ j_{r}}&& (]X[_{V^{r+1}})_{\an} \ar[d]^{p_{V^{r+1}}}\\
(X_V/C,O)_{\AN^\dagger} \ar[rr]^{{p_{X_V/C,O}} } && (]C[_{O})_{\an}
}
$$
from which we get a canonical isomorphism
$$
p_{X_V/C,O*}  j_{r*}j_{r}^{-1}\mathcal F =  p_{V^{r+1}*} \varphi_{V^{r+1}*}j_{r}^{-1}\mathcal F = p_{V^{r+1}*} \mathcal F_{V^{r+1}} = p_{V*}p_{r*} \mathcal F_{V^{r+1}}.
$$
giving a natural isomorphism of complexes
$$
p_{V/C,O*} \mathcal C(\mathcal F) \simeq p_{V*}\mathcal {CA}(\mathcal F).
$$

In order to prove the second assertion, we may choose an injective resolution $\mathcal I$ of $\mathcal F$.
By definition, $R\mathcal {CA}(\mathcal F) = [\mathcal {CA}(\mathcal I)]$ and all the terms of the latter are injective.
Thus, we have in the derived category
$$
Rp_{X_V/C,O*} \mathcal F  =  p_{X_V/C,O*} \mathcal I = p_{V*}[\mathcal {CA}(\mathcal I)] 
$$
$$
= Rp_{V*}[\mathcal {CA}(\mathcal I)] =  Rp_{V*} R\mathcal {CA}(\mathcal F). \quad \Box
$$

\begin{cor}
Let $S$ be a formal $\mathcal V$-scheme and $X \hookrightarrow P$ a good admissible formal embedding over $S$ with $P$ proper and smooth at $X$ over $S$.
If $\mathcal F$ is a complex of abelian sheaves on $X/S$ and $p : P \to S$ denotes the structural map, we have an isomorphism
$$
Rp_{X/S*} \mathcal F = Rp_{K*} R\mathcal {CA}(\mathcal F).
$$
\end{cor}

\textbf {Proof : }
We apply the proposition to the case $(C, O) = (S_{k}, S_{K}) =: S$ and $(X, V) = (X, P_{K}) =: (X \subset P)$.
Then, we saw in corollary 9.10 of \cite{LeStum06*} that $(X_{P}/S)_{\AN^\dagger} = (X/S)_{\AN^\dagger}$.
$\Box$

Of course, there is an analog of the proposition and the corollary for relative cohomology.

%
%

%
%
\bibliographystyle{plain}
\bibliography{BiblioBLS}

\begin{thebibliography}{1}

\bibitem{Berkovich93}
Vladimir~G. Berkovich.
\newblock \'{E}tale cohomology for non-{A}rchimedean analytic spaces.
\newblock {\em Inst. Hautes \'Etudes Sci. Publ. Math.}, (78):5--161 (1994),
  1993.

\bibitem{Berkovich99}
Vladimir~G. Berkovich.
\newblock Smooth {$p$}-adic analytic spaces are locally contractible.
\newblock {\em Invent. Math.}, 137(1):1--84, 1999.

\bibitem{Berthelot96c*}
Pierre Berthelot.
\newblock Cohomologie rigide et cohomologie \`a support propre, premi\`ere
  partie.
\newblock {\em Pr\'epublication de l'IRMAR}, 96--03:1--89, 1996.

\bibitem{Grothendieck68}
A.~Grothendieck.
\newblock Crystals and the de {R}ham cohomology of schemes.
\newblock In {\em Dix Expos\'es sur la Cohomologie des Sch\'emas}, pages
  306--358. North-Holland, Amsterdam, 1968.

\bibitem{KashiwaraSchapira90}
Masaki Kashiwara and Pierre Schapira.
\newblock {\em Sheaves on manifolds}, volume 292 of {\em Grundlehren der
  Mathematischen Wissenschaften [Fundamental Principles of Mathematical
  Sciences]}.
\newblock Springer-Verlag, Berlin, 1990.
\newblock With a chapter in French by Christian Houzel.

\bibitem{LeStum06*}
Bernard Le~Stum.
\newblock The overconvergent site 1. coefficients.
\newblock {\em Pr\'epublication de l'IRMAR}, 06(28):53, 2006.

\bibitem{LeStum07*}
Bernard Le~Stum.
\newblock {\em Rigid Cohomology}.
\newblock Cambridge Tracts in Mathematics. Cambridge University Press, 2007.

\bibitem{GrusonRaynaud71}
Michel Raynaud and Laurent Gruson.
\newblock Crit\`eres de platitude et de projectivit\'e. {T}echniques de
  ``platification'' d'un module.
\newblock {\em Invent. Math.}, 13:1--89, 1971.

\end{thebibliography}
%
%
\end{document}